\documentclass[11pt]{amsart}
\usepackage{amsmath,amssymb}
\usepackage[utf8,utf8x]{inputenc}
\usepackage[pagebackref=true,
pdftitle = {Configurations\ of\ lines},
pdfkeywords = {arrangement\ of\ lines;configuration\ of\ lines},
pdfsubject = {52C30}
pdfauthor = {}]
{hyperref}
\usepackage{hyperref}
\usepackage{longtable}

\usepackage{graphicx}
\usepackage{pict2e}

\newtheorem{propo}{Proposition}[section]

\newtheorem{algor}[propo]{Algorithm}

\theoremstyle{definition}

\theoremstyle{remark}
\newtheorem{remar}[propo]{Remark}

\numberwithin{equation}{section}

\newcommand{\algo}[6]
{
\begin{algor}{{\tt #1}{\rm (#2)}}\label{#1}\end{algor}
\vspace{-6pt}\noindent{\it #3}.

{\bf Input:} #4

{\bf Output:} #5

\newcounter{#1}
\begin{list}{\textbf{\arabic{#1}.}}{\usecounter{#1}}
#6\end{list}\vspace{3pt}}

\newcommand{\NN }{\mathbb{N}}
\newcommand{\CC }{\mathbb{C}}
\newcommand{\RR }{\mathbb{R}}
\newcommand{\FF }{\mathbb{F}}
\newcommand{\QQ }{\mathbb{Q}}
\newcommand{\ZZ }{\mathbb{Z}}
\newcommand{\PP }{\mathbb{P}}
\newcommand{\id }{\mathrm{id}}
\newcommand{\ii }{\mathrm{i}}
\newcommand{\Pc }{\mathcal{P}}
\newcommand{\Rc }{\mathcal{R}}
\newcommand{\Ac }{\mathcal{A}}
\newcommand{\Ec }{\mathcal{E}}
\newcommand{\Bc }{\mathcal{B}}
\newcommand{\Dc }{\mathcal{D}}
\newcommand{\Kc }{\mathcal{K}}
\newcommand{\Vc }{\mathcal{V}}
\newcommand{\Wc }{\mathcal{W}}
\newcommand{\Cc }{\mathcal{C}}
\newcommand{\Ie }{\mathfrak I}
\DeclareMathOperator{\diag}{diag}
\DeclareMathOperator{\lcm}{lcm}
\DeclareMathOperator{\RRe}{Re}
\DeclareMathOperator{\IIm}{Im}
\DeclareMathOperator{\Aut}{Aut}
\DeclareMathOperator{\End}{End}
\DeclareMathOperator{\Hom}{Hom}
\DeclareMathOperator{\GL}{GL}
\DeclareMathOperator{\PGL}{PGL}
\DeclareMathOperator{\PGAL}{P\Gamma L}
\DeclareMathOperator{\Gal}{Gal}
\DeclareMathOperator{\PG}{PG}
\DeclareMathOperator{\rank}{rank}
\DeclareMathOperator{\Alt}{Alt}
\DeclareMathOperator{\md}{mod}
\DeclareMathOperator{\mm}{m}
\newcommand{\Spairs }{\Psi}
\newcommand{\Sc }{\mathcal{S}}
\newcommand{\Sym }{\mathbb{S}}
\newcommand{\PFq }{\FF_q\PP^2}

\newcommand{\IF }{\mathcal{IF}}
\newcommand{\RF }{\mathcal{RF}}
\DeclareMathOperator{\Der}{Der}
\DeclareMathOperator{\pdeg}{pdeg}
\DeclareMathOperator{\rk}{rk}

\newcommand{\lmultiset }{\{\hspace{-4.7pt}\{}
\newcommand{\rmultiset }{\}\hspace{-4.7pt}\}}

\title[$(22_4)$ and $(26_4)$ configurations of lines]
{$(22_4)$ and $(26_4)$ configurations of lines}

\author{M.~Cuntz}
\address{Michael Cuntz,
Institut f\"ur Algebra, Zahlentheorie und Diskrete Mathematik,
Fakult\"at f\"ur Mathematik und Physik,
Leibniz Universit\"at Hannover,
Welfengarten 1,
D-30167 Hannover, Germany}
\email{cuntz@math.uni-hannover.de}

\begin{document}

\begin{abstract}
We present a technique to produce arrangements of lines with nice properties.
As an application, we construct $(22_4)$ and $(26_4)$ configurations of lines.
Thus concerning the existence of geometric $(n_4)$ configurations, only the case $n=23$ remains open.
\end{abstract}


\maketitle

\section{Enumerating arrangements}

There are several ways to enumerate arrangements of lines in the real plane. For instance, one can enumerate all wiring diagrams and thus oriented matroids. However, without a very strong local condition on the cell structure, such an enumeration is feasible only for a small number of lines. In any case, most types of interesting arrangements of more than say $20$ lines can probably not be enumerated completely (nowadays by a computer).

A much more promising method is (as already noted by many authors) to exploit symmetry. In fact, most relevant examples in the literature have a non-trivial symmetry group. Symmetry reduces the degrees of freedom considerably and allows us to compute examples with many more lines.
The following (very simple) algorithm is a useful tool to produce ``interesting'' examples of arrangements with non-trivial symmetry group:
\algo{Enumerate arrangements}{$q$,$P$}{Look for matroids with $P$ which are realizable over $\CC$}
{a prime power $q$, a property $P$}{matroids of arrangements of lines in $\CC\PP^2$ with $P$}
{
\item Depending on $P$, choose a small set of lines $\Ac_0\subseteq \PFq$ and an $n\in\NN$.
\item For every group $H\le \PGL_3(\FF_q)$ with $|H|=n$, compute the orbit $\Ac:=H\Ac_0$.
\item If $\Ac$ has property $P$, then compute its matroid $M$.
Print $M$ if it is realizable over $\CC$.
}

\begin{remar}
\begin{enumerate}
\item If $q$ is not too big, then it is indeed possible to compute all the subgroups $H$ with $|H|=n$.
However, if $q$ is too small, then only very few matroids $M$ will be realizable in characteristic zero.
\item If we are looking for arrangements with $m = nk$ lines, then it is good to choose $\Ac_0$ with approximately $k$ lines.
\item This algorithm mostly produces matroids that are not orientable. Thus it is a priori not the best method if one is searching for arrangements in the real projective plane. On the other hand, most ``interesting'' arrangements will define a matroid that is realizable over many finite fields, such that these matroids will certainly appear in the enumeration.
\item Realizing rank three matroids with a small number of lines, depending on the matroid maybe up to 70 lines, is not easy but works in most cases (see for example \cite{p-C10b}).
\end{enumerate}
\end{remar}

\begin{figure}
\begin{center}
\setlength{\unitlength}{0.6pt}
\begin{picture}(400,400)(100,200)
\put(486.356823840916435087096077209,377.058426612943823409279042190){\circle{5}}
\put(352.324599306993152919454258755,393.558547467119200261987656865){\circle{5}}
\put(247.675400693006847080545741245,406.441452532880799738012343135){\circle{5}}
\put(113.643176159083564912903922791,422.941573387056176590720957810){\circle{5}}
\put(322.941573387056176590720957810,341.233948160126047342533398111){\circle{5}}
\put(218.292374773069870751812440301,318.292374773069870751812440301){\circle{5}}
\put(277.058426612943823409279042190,331.175279838831470227837126571){\circle{5}}
\put(358.766051839873952657466601890,377.058426612943823409279042190){\circle{5}}
\put(322.941573387056176590720957810,468.824720161168529772162873429){\circle{5}}
\put(277.058426612943823409279042190,586.356823840916435087096077209){\circle{5}}
\put(381.707625226930129248187559700,481.707625226930129248187559700){\circle{5}}
\put(277.058426612943823409279042190,422.941573387056176590720957810){\circle{5}}
\put(322.941573387056176590720957810,377.058426612943823409279042190){\circle{5}}
\put(231.175279838831470227837126571,377.058426612943823409279042190){\circle{5}}
\put(322.941573387056176590720957810,213.643176159083564912903922792){\circle{5}}
\put(368.824720161168529772162873429,422.941573387056176590720957810){\circle{5}}
\put(306.441452532880799738012343135,347.675400693006847080545741245){\circle{5}}
\put(241.233948160126047342533398110,422.941573387056176590720957810){\circle{5}}
\put(293.558547467119200261987656865,452.324599306993152919454258755){\circle{5}}
\put(277.058426612943823409279042190,458.766051839873952657466601890){\circle{5}}
\Line(498.501511336580603580962197971,375.563347260859106800755213289)(101.498488663419396419037802029,424.436652739140893199244786711)
\Line(498.979892613182403738163682238,379.825701111404176488213105521)(127.708402582995263278924224061,298.434230867399193737967475333)
\Line(413.970105744024693682231873460,235.650327056316204169928181623)(272.475233989393202875093764441,598.096913797417209490458480924)
\Line(456.951253954483618827639063207,523.960896584831098110927352719)(112.459981529441479051399659714,330.510853566455654721197687080)
\Line(498.679853559756570979912758077,377.058426612943823409279042190)(101.320146440243429020087241923,377.058426612943823409279042190)
\Line(322.941573387056176590720957810,598.679853559756570979912758077)(322.941573387056176590720957810,201.320146440243429020087241923)
\Line(437.499994154186245570075627493,545.236881017180853995327706099)(247.560520513175409383010689386,206.997147712395444582344768462)
\Line(487.540018470558520948600340286,469.489146433544345278802312920)(143.048746045516381172360936793,276.039103415168901889072647281)
\Line(493.002852287604555417655231538,452.439479486824590616989310613)(154.763118982819146004672293901,262.500005845813754429924372507)
\Line(445.236881017180853995327706099,537.499994154186245570075627493)(106.997147712395444582344768462,347.560520513175409383010689386)
\Line(327.524766010606797124906235560,201.903086202582790509541519075)(186.029894255975306317768126540,564.349672943683795830071818377)
\Line(401.565769132600806262032524667,572.291597417004736721075775939)(320.174298888595823511786894480,201.020107386817596261836317763)
\Line(279.825701111404176488213105521,598.979892613182403738163682238)(198.434230867399193737967475333,227.708402582995263278924224061)
\Line(472.291597417004736721075775939,501.565769132600806262032524667)(101.020107386817596261836317763,420.174298888595823511786894480)
\Line(423.960896584831098110927352719,556.951253954483618827639063207)(230.510853566455654721197687080,212.459981529441479051399659714)
\Line(498.096913797417209490458480924,372.475233989393202875093764441)(135.650327056316204169928181623,513.970105744024693682231873460)
\Line(498.679853559756570979912758077,422.941573387056176590720957810)(101.320146440243429020087241923,422.941573387056176590720957810)
\Line(277.058426612943823409279042190,598.679853559756570979912758077)(277.058426612943823409279042190,201.320146440243429020087241923)
\Line(324.436652739140893199244786711,201.498488663419396419037802029)(275.563347260859106800755213289,598.501511336580603580962197971)
\Line(352.439479486824590616989310613,593.002852287604555417655231538)(162.500005845813754429924372507,254.763118982819146004672293901)
\Line(464.349672943683795830071818377,286.029894255975306317768126539)(101.903086202582790509541519075,427.524766010606797124906235560)
\Line(369.489146433544345278802312920,587.540018470558520948600340285)(176.039103415168901889072647281,243.048746045516381172360936793)
\thinlines
\strokepath
\end{picture}
\vspace{20pt}

\begin{picture}(400,400)(100,200)
\put(359.561670289725574667590751435,380.146109903424808444136416189){\circle{5}}
\put(461.275246333884696436032456233,380.146109903424808444136416189){\circle{5}}
\put(229.289321881345247559915563789,380.146109903424808444136416189){\circle{5}}
\put(345.282284107614971997974010011,329.289321881345247559915563789){\circle{5}}
\put(410.418458311805135551811603833,329.289321881345247559915563789){\circle{5}}
\put(280.146109903424808444136416189,329.289321881345247559915563789){\circle{5}}
\put(319.853890096575191555863583811,470.710678118654752440084436211){\circle{5}}
\put(189.581541688194864448188396167,470.710678118654752440084436211){\circle{5}}
\put(254.717715892385028002025989989,470.710678118654752440084436211){\circle{5}}
\put(370.710678118654752440084436211,419.853890096575191555863583811){\circle{5}}
\put(240.438329710274425332409248566,419.853890096575191555863583811){\circle{5}}
\put(138.724753666115303563967543768,419.853890096575191555863583811){\circle{5}}
\put(319.853890096575191555863583811,238.724753666115303563967543768){\circle{5}}
\put(280.146109903424808444136416189,459.561670289725574667590751435){\circle{5}}
\put(229.289321881345247559915563789,510.418458311805135551811603833){\circle{5}}
\put(370.710678118654752440084436211,354.717715892385028002025989989){\circle{5}}
\put(280.146109903424808444136416189,561.275246333884696436032456233){\circle{5}}
\put(319.853890096575191555863583811,340.438329710274425332409248566){\circle{5}}
\put(229.289321881345247559915563789,445.282284107614971997974010011){\circle{5}}
\put(370.710678118654752440084436211,289.581541688194864448188396167){\circle{5}}
\Line(499.012117842188477954921988695,380.146109903424808444136416189)(100.987882157811522045078011305,380.146109903424808444136416189)
\Line(487.082869338697069279187436616,329.289321881345247559915563789)(112.917130661302930720812563384,329.289321881345247559915563789)
\Line(487.082869338697069279187436616,470.710678118654752440084436211)(112.917130661302930720812563384,470.710678118654752440084436211)
\Line(499.012117842188477954921988695,419.853890096575191555863583811)(100.987882157811522045078011305,419.853890096575191555863583811)
\Line(411.714493677880360916003934651,565.891144737429637376841996797)(309.038263942426231078593425378,200.204329914517322033677434938)
\Line(459.874685152285089959924033189,279.833095040865293151803134434)(179.833095040865293151803134434,559.874685152285089959924033188)
\Line(420.757884611552781857849372664,240.571228116933348049012040785)(264.555410330033193795717058860,596.834146080215671435389586479)
\Line(499.795670085482677966322565062,390.961736057573768921406574622)(134.108855262570362623158003204,288.285506322119639083996065349)
\Line(335.444589669966806204282941140,203.165853919784328564610413521)(179.242115388447218142150627336,559.428771883066651950987959215)
\Line(280.146109903424808444136416189,599.012117842188477954921988695)(280.146109903424808444136416189,200.987882157811522045078011305)
\Line(468.424587950951963797615435626,292.142880747017693320801326074)(131.575412049048036202384564374,507.857119252982306679198673925)
\Line(420.166904959134706848196865566,240.125314847714910040075966811)(140.125314847714910040075966811,520.166904959134706848196865566)
\Line(318.054302078973365566221202889,599.183438509433253558117242933)(100.816561490566746441882757067,381.945697921026634433778797111)
\Line(499.183438509433253558117242933,418.054302078973365566221202889)(281.945697921026634433778797111,200.816561490566746441882757067)
\Line(496.834146080215671435389586479,364.555410330033193795717058860)(140.571228116933348049012040785,520.757884611552781857849372664)
\Line(319.853890096575191555863583811,599.012117842188477954921988695)(319.853890096575191555863583811,200.987882157811522045078011305)
\Line(459.428771883066651950987959215,279.242115388447218142150627336)(103.165853919784328564610413521,435.444589669966806204282941140)
\Line(407.857119252982306679198673925,231.575412049048036202384564374)(192.142880747017693320801326074,568.424587950951963797615435627)
\Line(370.710678118654752440084436211,587.082869338697069279187436616)(370.710678118654752440084436211,212.917130661302930720812563384)
\Line(465.891144737429637376841996797,511.714493677880360916003934651)(100.204329914517322033677434938,409.038263942426231078593425378)
\Line(229.289321881345247559915563789,587.082869338697069279187436616)(229.289321881345247559915563789,212.917130661302930720812563384)
\Line(290.961736057573768921406574622,599.795670085482677966322565062)(188.285506322119639083996065349,234.108855262570362623158003204)
\thinlines
\strokepath
\end{picture}
\end{center}
\caption{Two dual $22_4$-configurations of lines ($w=\frac{-7+3\sqrt{17}}{2}$).\label{fig_22_4}}
\end{figure}
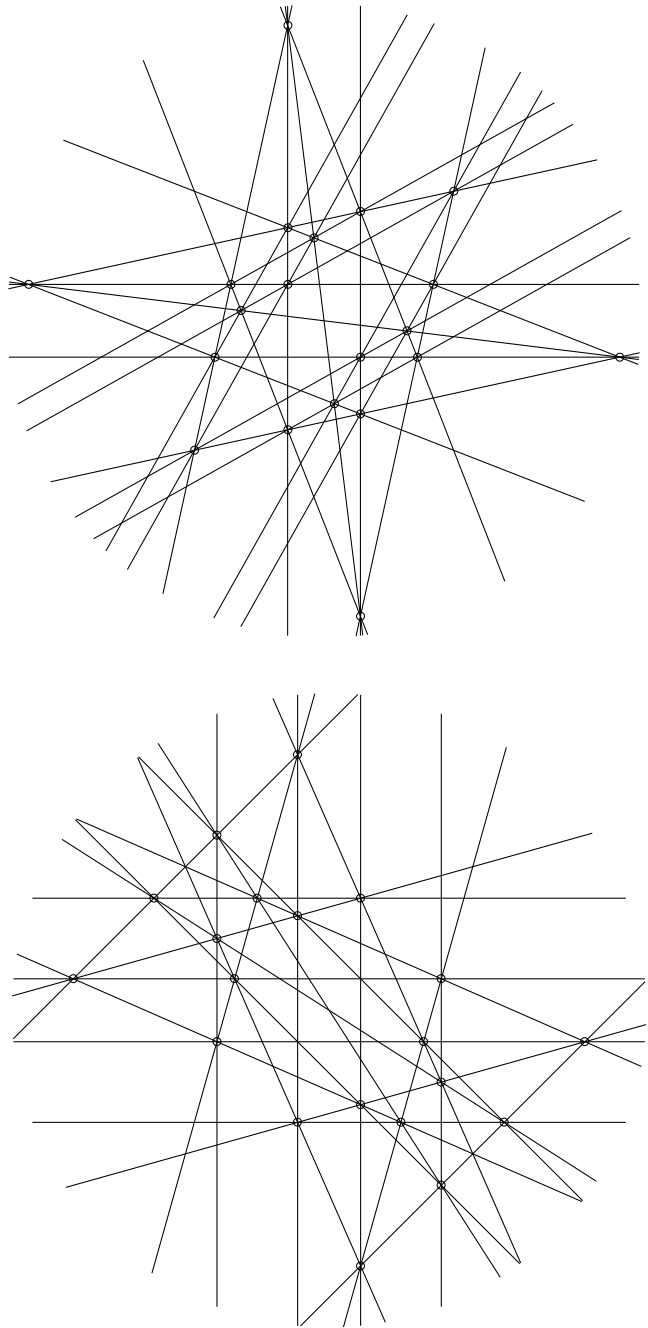

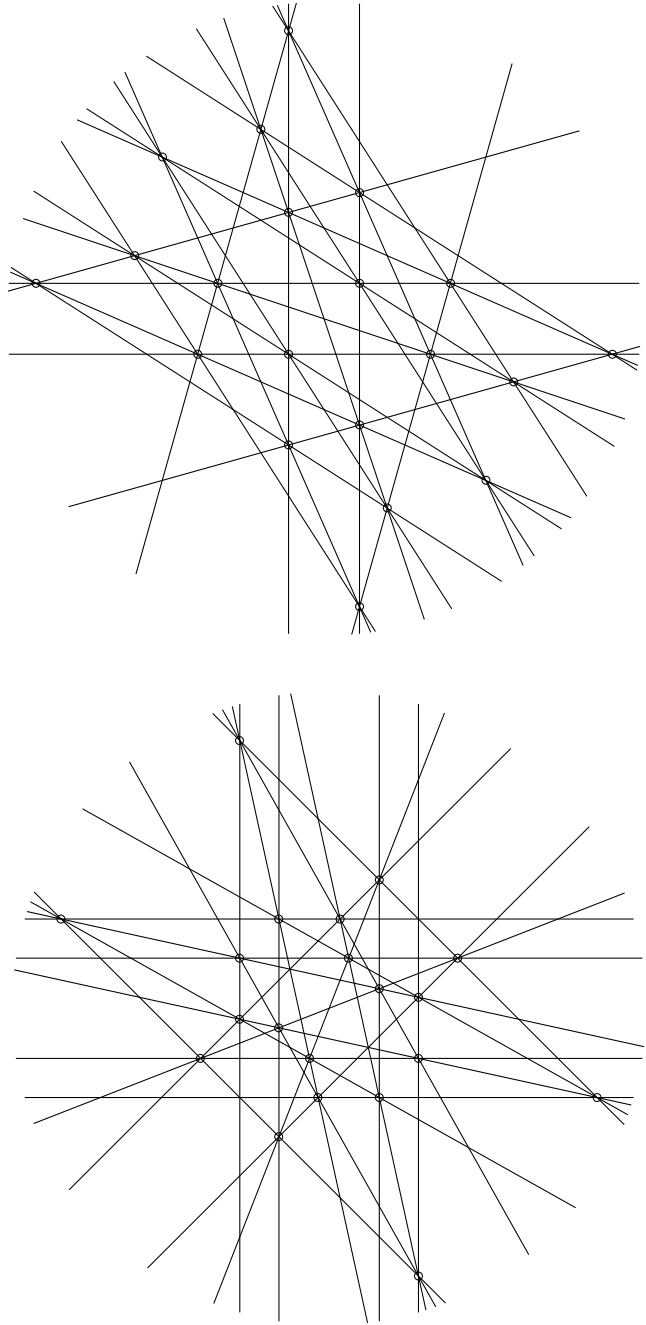
\begin{figure}
\begin{center}
\setlength{\unitlength}{0.6pt}
\begin{picture}(400,400)(100,200)
\put(367.082039324993690892275210062,377.639320225002103035908263313){\circle{5}}
\put(180.541887076558077005776535340,439.819370974480640998074488220){\circle{5}}
\put(419.458112923441922994223464660,360.180629025519359001925511780){\circle{5}}
\put(232.917960675006309107724789938,422.360679774997896964091736687){\circle{5}}
\put(220.361258051038718003851023560,377.639320225002103035908263313){\circle{5}}
\put(277.639320225002103035908263313,377.639320225002103035908263313){\circle{5}}
\put(481.638163672920460956389689567,377.639320225002103035908263313){\circle{5}}
\put(322.360679774997896964091736687,479.638741948961281996148976440){\circle{5}}
\put(118.361836327079539043610310433,422.360679774997896964091736687){\circle{5}}
\put(277.639320225002103035908263313,467.082039324993690892275210062){\circle{5}}
\put(322.360679774997896964091736687,422.360679774997896964091736687){\circle{5}}
\put(198.000578276040821039759286873,501.999421723959178960240713127){\circle{5}}
\put(401.999421723959178960240713127,298.000578276040821039759286873){\circle{5}}
\put(277.639320225002103035908263313,581.638163672920460956389689567){\circle{5}}
\put(322.360679774997896964091736687,332.917960675006309107724789938){\circle{5}}
\put(322.360679774997896964091736687,218.361836327079539043610310433){\circle{5}}
\put(260.180629025519359001925511780,519.458112923441922994223464660){\circle{5}}
\put(277.639320225002103035908263313,320.361258051038718003851023560){\circle{5}}
\put(339.819370974480640998074488220,280.541887076558077005776535340){\circle{5}}
\put(379.638741948961281996148976440,422.360679774997896964091736687){\circle{5}}
\Line(489.736659610102759919933612666,336.754446796632413360022129111)(110.263340389897240080066387334,463.245553203367586639977870889)
\Line(498.746069143517905200435966577,377.639320225002103035908263313)(101.253930856482094799564033423,377.639320225002103035908263313)
\Line(461.054672240373943871890584238,518.580742743270597912524823203)(100.759952976983799440071255927,417.418486222008307367796127633)
\Line(483.061814743914600319592527135,319.449568699697882651735716776)(150.254423113298622607231348655,532.575496238441461918269182668)
\Line(425.540275315678140229647971569,244.309154817427579652556050738)(270.544837448586885134247112298,597.819092605010779433045861937)
\Line(455.690845182572420347443949262,274.459724684321859770352028431)(102.180907394989220566954138063,429.455162551413114865752887702)
\Line(332.424950472861671010762344439,202.645946109960873452770295163)(134.310098098285437281942152403,512.012751094686559050400929605)
\Line(449.745576886701377392768651345,267.424503761558538081730817333)(116.938185256085399680407472865,480.550431300302117348264283224)
\Line(497.354053890039126547229704837,367.575049527138328989237655561)(187.987248905313440949599070395,565.689901901714562718057847597)
\Line(412.012751094686559050400929605,234.310098098285437281942152403)(102.645946109960873452770295163,432.424950472861671010762344439)
\Line(499.240047023016200559928744073,382.581513777991692632203872367)(138.945327759626056128109415762,281.419257256729402087475176797)
\Line(322.360679774997896964091736687,598.746069143517905200435966577)(322.360679774997896964091736687,201.253930856482094799564033423)
\Line(277.639320225002103035908263313,598.746069143517905200435966577)(277.639320225002103035908263313,201.253930856482094799564033423)
\Line(498.746069143517905200435966577,422.360679774997896964091736687)(101.253930856482094799564033423,422.360679774997896964091736687)
\Line(432.575496238441461918269182667,250.254423113298622607231348655)(219.449568699697882651735716776,583.061814743914600319592527134)
\Line(418.580742743270597912524823203,561.054672240373943871890584238)(317.418486222008307367796127633,200.759952976983799440071255927)
\Line(329.455162551413114865752887702,202.180907394989220566954138063)(174.459724684321859770352028431,555.690845182572420347443949262)
\Line(497.819092605010779433045861937,370.544837448586885134247112298)(144.309154817427579652556050738,525.540275315678140229647971569)
\Line(363.245553203367586639977870889,210.263340389897240080066387334)(236.754446796632413360022129111,589.736659610102759919933612666)
\Line(465.689901901714562718057847597,287.987248905313440949599070395)(267.575049527138328989237655561,597.354053890039126547229704837)
\Line(282.581513777991692632203872367,599.240047023016200559928744073)(181.419257256729402087475176797,238.945327759626056128109415762)
\Line(380.550431300302117348264283224,216.938185256085399680407472865)(167.424503761558538081730817333,549.745576886701377392768651345)
\thinlines
\strokepath
\end{picture}
\vspace{20pt}

\begin{picture}(400,400)(100,200)
\put(131.060716568171734311280321889,456.313094477276088562906559370){\circle{5}}
\put(306.932458726091498077071311518,456.313094477276088562906559370){\circle{5}}
\put(268.377223398316206680011064556,456.313094477276088562906559370){\circle{5}}
\put(287.654841062203852378541188037,368.377223398316206680011064556){\circle{5}}
\put(218.996587647131616194175816704,368.377223398316206680011064556){\circle{5}}
\put(356.313094477276088562906559370,368.377223398316206680011064556){\circle{5}}
\put(243.686905522723911437093440630,431.622776601683793319988935444){\circle{5}}
\put(381.003412352868383805824183296,431.622776601683793319988935444){\circle{5}}
\put(312.345158937796147621458811963,431.622776601683793319988935444){\circle{5}}
\put(331.622776601683793319988935444,343.686905522723911437093440630){\circle{5}}
\put(468.939283431828265688719678111,343.686905522723911437093440630){\circle{5}}
\put(293.067541273908501922928688481,343.686905522723911437093440630){\circle{5}}
\put(243.686905522723911437093440630,393.067541273908501922928688481){\circle{5}}
\put(356.313094477276088562906559370,231.060716568171734311280321889){\circle{5}}
\put(268.377223398316206680011064556,318.996587647131616194175816704){\circle{5}}
\put(331.622776601683793319988935444,412.345158937796147621458811963){\circle{5}}
\put(356.313094477276088562906559370,406.932458726091498077071311518){\circle{5}}
\put(243.686905522723911437093440630,568.939283431828265688719678111){\circle{5}}
\put(268.377223398316206680011064556,387.654841062203852378541188037){\circle{5}}
\put(331.622776601683793319988935444,481.003412352868383805824183296){\circle{5}}
\Line(491.908403647139374806395530474,456.313094477276088562906559370)(108.091596352860625193604469526,456.313094477276088562906559370)
\Line(497.484176581314990174384610437,368.377223398316206680011064556)(102.515823418685009825615389563,368.377223398316206680011064556)
\Line(497.484176581314990174384610437,431.622776601683793319988935444)(102.515823418685009825615389563,431.622776601683793319988935444)
\Line(491.908403647139374806395530474,343.686905522723911437093440630)(108.091596352860625193604469526,343.686905522723911437093440630)
\Line(455.444300125756513269905714769,274.154580701506128944891118697)(111.635276040120426406204119997,467.221505245853949799845218093)
\Line(373.412906734483218844461093281,213.960904310964604029725787978)(113.960904310964604029725787978,473.412906734483218844461093281)
\Line(498.525159788499978712131880586,375.756218715914980280545111807)(109.539479100724373218822407106,461.031057493513842858803943657)
\Line(425.845419298493871055108881303,244.555699874243486730094285232)(232.778494754146050200154781907,588.364723959879573593795880003)
\Line(490.460520899275626781177592894,338.968942506486157141196056343)(101.474840211500021287868119414,424.243781284085019719454888193)
\Line(356.313094477276088562906559370,591.908403647139374806395530475)(356.313094477276088562906559370,208.091596352860625193604469526)
\Line(486.306417630942395563794016946,472.731827623983686748180653471)(113.693582369057604436205983054,327.268172376016313251819346529)
\Line(486.039095689035395970274212022,326.587093265516781155538906719)(226.587093265516781155538906719,586.039095689035395970274212022)
\Line(463.939692032901194685372519229,514.559056281716604199537271376)(185.440943718283395800462728624,236.060307967098805314627480771)
\Line(414.559056281716604199537271376,563.939692032901194685372519229)(136.060307967098805314627480771,285.440943718283395800462728624)
\Line(361.031057493513842858803943657,209.539479100724373218822407106)(275.756218715914980280545111807,598.525159788499978712131880586)
\Line(243.686905522723911437093440630,591.908403647139374806395530475)(243.686905522723911437093440630,208.091596352860625193604469526)
\Line(324.243781284085019719454888193,201.474840211500021287868119414)(238.968942506486157141196056343,590.460520899275626781177592894)
\Line(372.731827623983686748180653471,586.306417630942395563794016945)(227.268172376016313251819346529,213.693582369057604436205983054)
\Line(331.622776601683793319988935444,597.484176581314990174384610437)(331.622776601683793319988935444,202.515823418685009825615389563)
\Line(367.221505245853949799845218093,211.635276040120426406204119997)(174.154580701506128944891118697,555.444300125756513269905714769)
\Line(268.377223398316206680011064556,597.484176581314990174384610437)(268.377223398316206680011064556,202.515823418685009825615389563)
\Line(488.364723959879573593795880003,332.778494754146050200154781906)(144.555699874243486730094285232,525.845419298493871055108881304)
\thinlines
\strokepath
\end{picture}
\end{center}
\caption{Two dual $22_4$-configurations of lines ($w=\frac{-7-3\sqrt{17}}{2}$).\label{fig_22_4_e}}
\end{figure}

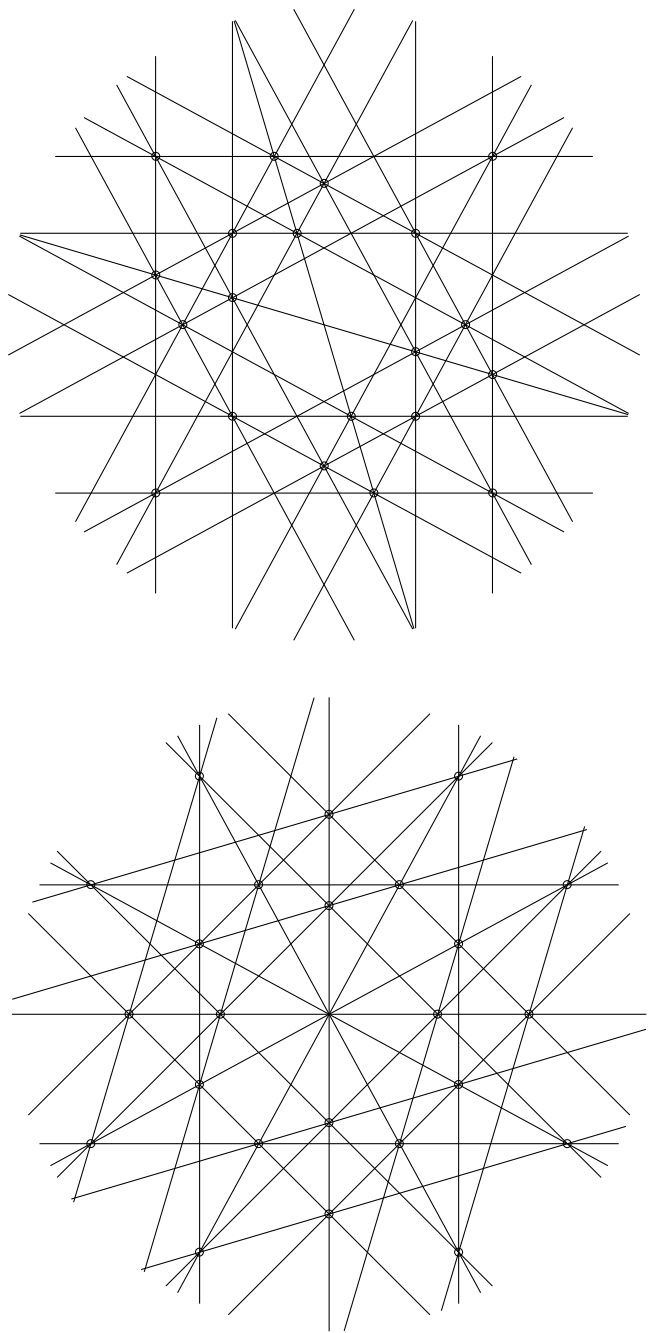
\begin{figure}
\begin{center}
\setlength{\unitlength}{0.6pt}
\begin{picture}(400,400)(100,200)
\put(406.191270323980923959719363144,368.610100216678785379156576753){\circle{5}}
\put(193.808729676019076040280636856,431.389899783321214620843423248){\circle{5}}
\put(242.264973081037423549085121950,417.066343621697132887961061847){\circle{5}}
\put(357.735026918962576450914878050,382.933656378302867112038938153){\circle{5}}
\put(317.066343621697132887961061847,342.264973081037423549085121950){\circle{5}}
\put(282.933656378302867112038938153,457.735026918962576450914878050){\circle{5}}
\put(331.389899783321214620843423248,293.808729676019076040280636856){\circle{5}}
\put(268.610100216678785379156576753,506.191270323980923959719363144){\circle{5}}
\put(242.264973081037423549085121950,342.264973081037423549085121950){\circle{5}}
\put(357.735026918962576450914878050,342.264973081037423549085121950){\circle{5}}
\put(357.735026918962576450914878050,457.735026918962576450914878050){\circle{5}}
\put(389.124926702283791071758301298,400.000000000000000000000000000){\circle{5}}
\put(210.875073297716208928241698702,400.000000000000000000000000000){\circle{5}}
\put(193.808729676019076040280636856,506.191270323980923959719363144){\circle{5}}
\put(193.808729676019076040280636856,293.808729676019076040280636856){\circle{5}}
\put(406.191270323980923959719363144,293.808729676019076040280636856){\circle{5}}
\put(406.191270323980923959719363144,506.191270323980923959719363144){\circle{5}}
\put(242.264973081037423549085121950,457.735026918962576450914878050){\circle{5}}
\put(300.000000000000000000000000000,310.875073297716208928241698702){\circle{5}}
\put(300.000000000000000000000000000,489.124926702283791071758301298){\circle{5}}
\Line(491.796076874067532766006639931,343.305512651433408089006258135)(108.203923125932467233993360069,456.694487348566591910993741865)
\Line(356.694487348566591910993741865,208.203923125932467233993360069)(243.305512651433408089006258135,591.796076874067532766006639931)
\Line(491.485421551267621995020382274,342.264973081037423549085121950)(108.514578448732378004979617726,342.264973081037423549085121950)
\Line(456.696264548597265936791884181,275.716933267174245861507966481)(280.884905558704872643770902211,599.084437273485463477367973416)
\Line(319.115094441295127356229097789,200.915562726514536522632026584)(143.303735451402734063208115819,524.283066732825754138492033519)
\Line(193.808729676019076040280636856,569.479833924214146761498507626)(193.808729676019076040280636856,230.520166075785853238501492374)
\Line(406.191270323980923959719363144,569.479833924214146761498507626)(406.191270323980923959719363144,230.520166075785853238501492374)
\Line(242.264973081037423549085121950,591.485421551267621995020382274)(242.264973081037423549085121950,208.514578448732378004979617726)
\Line(355.941343290630579660618079662,207.982901514370806724901045967)(169.257286168709711000505980441,551.348415188363749712145137829)
\Line(430.742713831290288999494019559,248.651584811636250287854862171)(244.058656709369420339381920338,592.017098485629193275098954033)
\Line(357.735026918962576450914878050,591.485421551267621995020382274)(357.735026918962576450914878050,208.514578448732378004979617726)
\Line(492.017098485629193275098954033,344.058656709369420339381920338)(148.651584811636250287854862171,530.742713831290288999494019559)
\Line(492.017098485629193275098954033,455.941343290630579660618079662)(148.651584811636250287854862171,269.257286168709711000505980441)
\Line(451.348415188363749712145137829,269.257286168709711000505980441)(107.982901514370806724901045967,455.941343290630579660618079662)
\Line(499.084437273485463477367973415,419.115094441295127356229097789)(175.716933267174245861507966481,243.303735451402734063208115819)
\Line(491.485421551267621995020382274,457.735026918962576450914878050)(108.514578448732378004979617726,457.735026918962576450914878050)
\Line(456.696264548597265936791884181,524.283066732825754138492033519)(280.884905558704872643770902211,200.915562726514536522632026584)
\Line(430.742713831290288999494019559,551.348415188363749712145137829)(244.058656709369420339381920338,207.982901514370806724901045967)
\Line(469.479833924214146761498507626,293.808729676019076040280636856)(130.520166075785853238501492374,293.808729676019076040280636856)
\Line(355.941343290630579660618079662,592.017098485629193275098954033)(169.257286168709711000505980441,248.651584811636250287854862171)
\Line(451.348415188363749712145137829,530.742713831290288999494019559)(107.982901514370806724901045967,344.058656709369420339381920338)
\Line(424.283066732825754138492033519,243.303735451402734063208115819)(100.915562726514536522632026584,419.115094441295127356229097789)
\Line(499.084437273485463477367973415,380.884905558704872643770902211)(175.716933267174245861507966481,556.696264548597265936791884181)
\Line(424.283066732825754138492033519,556.696264548597265936791884181)(100.915562726514536522632026584,380.884905558704872643770902211)
\Line(469.479833924214146761498507626,506.191270323980923959719363144)(130.520166075785853238501492374,506.191270323980923959719363144)
\Line(319.115094441295127356229097789,599.084437273485463477367973416)(143.303735451402734063208115819,275.716933267174245861507966481)
\thinlines
\strokepath
\end{picture}
\vspace{20pt}

\begin{picture}(400,400)(100,200)
\put(381.649658092772603273242802490,355.607978004894860689352196086){\circle{5}}
\put(255.607978004894860689352196086,318.350341907227396726757197510){\circle{5}}
\put(300.000000000000000000000000000,331.472523394971209596277955014){\circle{5}}
\put(344.392021995105139310647803914,481.649658092772603273242802490){\circle{5}}
\put(300.000000000000000000000000000,468.527476605028790403722044986){\circle{5}}
\put(218.350341907227396726757197510,444.392021995105139310647803914){\circle{5}}
\put(300.000000000000000000000000000,526.041680087877742583890606404){\circle{5}}
\put(381.649658092772603273242802490,550.177134697801393676964847476){\circle{5}}
\put(149.822865302198606323035152524,481.649658092772603273242802490){\circle{5}}
\put(300.000000000000000000000000000,273.958319912122257416109393596){\circle{5}}
\put(218.350341907227396726757197510,249.822865302198606323035152524){\circle{5}}
\put(450.177134697801393676964847476,318.350341907227396726757197510){\circle{5}}
\put(426.041680087877742583890606404,400.000000000000000000000000000){\circle{5}}
\put(450.177134697801393676964847476,481.649658092772603273242802490){\circle{5}}
\put(381.649658092772603273242802490,249.822865302198606323035152524){\circle{5}}
\put(173.958319912122257416109393596,400.000000000000000000000000000){\circle{5}}
\put(149.822865302198606323035152524,318.350341907227396726757197510){\circle{5}}
\put(218.350341907227396726757197510,550.177134697801393676964847476){\circle{5}}
\put(381.649658092772603273242802490,444.392021995105139310647803914){\circle{5}}
\put(344.392021995105139310647803914,318.350341907227396726757197510){\circle{5}}
\put(368.527476605028790403722044986,400.000000000000000000000000000){\circle{5}}
\put(231.472523394971209596277955014,400.000000000000000000000000000){\circle{5}}
\put(218.350341907227396726757197510,355.607978004894860689352196086){\circle{5}}
\put(255.607978004894860689352196086,481.649658092772603273242802490){\circle{5}}
\Line(499.775468099116178528439127268,390.525700776362713675101904630)(137.482167998551285434155871308,283.432619135759543741007488966)
\Line(462.517832001448714565844128692,516.567380864240456258992511034)(100.224531900883821471560872732,409.474299223637286324898095370)
\Line(418.542739441042961225234296967,561.082646259033743019600284479)(112.929783953928248371043658046,470.744146531540253930607365488)
\Line(487.070216046071751628956341954,329.255853468459746069392634512)(181.457260558957038774765703033,238.917353740966256980399715521)
\Line(461.082646259033743019600284479,518.542739441042961225234296967)(370.744146531540253930607365488,212.929783953928248371043658046)
\Line(229.255853468459746069392634512,587.070216046071751628956341955)(138.917353740966256980399715521,281.457260558957038774765703033)
\Line(416.567380864240456258992511033,562.517832001448714565844128692)(309.474299223637286324898095370,200.224531900883821471560872732)
\Line(290.525700776362713675101904630,599.775468099116178528439127268)(183.432619135759543741007488967,237.482167998551285434155871308)
\Line(489.624053747842418146039029213,336.417626340035324437851577191)(236.417626340035324437851577191,589.624053747842418146039029213)
\Line(489.624053747842418146039029213,463.582373659964675562148422808)(236.417626340035324437851577191,210.375946252157581853960970788)
\Line(500.000000000000000000000000000,400.000000000000000000000000000)(100.000000000000000000000000000,400.000000000000000000000000000)
\Line(395.531250102563154769060803914,575.709363022695882833807827128)(204.468749897436845230939196086,224.290636977304117166192172871)
\Line(363.582373659964675562148422808,210.375946252157581853960970788)(110.375946252157581853960970788,463.582373659964675562148422808)
\Line(482.574185835055371152323260934,318.350341907227396726757197510)(117.425814164944628847676739066,318.350341907227396726757197510)
\Line(482.574185835055371152323260934,481.649658092772603273242802490)(117.425814164944628847676739066,481.649658092772603273242802490)
\Line(300.000000000000000000000000000,600.000000000000000000000000000)(300.000000000000000000000000000,200.000000000000000000000000000)
\Line(402.944119480984316310573587501,571.471596086013106714295632487)(128.528403913986893285704367513,297.055880519015683689426412499)
\Line(471.471596086013106714295632487,297.055880519015683689426412499)(197.055880519015683689426412499,571.471596086013106714295632487)
\Line(471.471596086013106714295632487,502.944119480984316310573587501)(197.055880519015683689426412499,228.528403913986893285704367513)
\Line(402.944119480984316310573587501,228.528403913986893285704367513)(128.528403913986893285704367513,502.944119480984316310573587501)
\Line(363.582373659964675562148422808,589.624053747842418146039029213)(110.375946252157581853960970788,336.417626340035324437851577191)
\Line(381.649658092772603273242802490,582.574185835055371152323260934)(381.649658092772603273242802490,217.425814164944628847676739066)
\Line(218.350341907227396726757197510,582.574185835055371152323260934)(218.350341907227396726757197510,217.425814164944628847676739066)
\Line(475.709363022695882833807827128,495.531250102563154769060803914)(124.290636977304117166192172871,304.468749897436845230939196086)
\Line(475.709363022695882833807827128,304.468749897436845230939196086)(124.290636977304117166192172871,495.531250102563154769060803914)
\Line(395.531250102563154769060803914,224.290636977304117166192172871)(204.468749897436845230939196086,575.709363022695882833807827128)
\thinlines
\strokepath
\end{picture}
\end{center}
\caption{Two dual $26_4$-configurations of lines.\label{fig_26_4}}
\end{figure}

\section{$(n_k)$ configurations of lines}

A configuration of lines and points is an \emph{$(n_k)$ configuration} if it consists of $n$ lines and $n$ points, each of which is incident to exactly $k$ of the other type. It is called \emph{geometric} if these are points and lines in the real projective plane.

There are many results concerning geometric $(n_4)$ configurations:
\begin{enumerate}
\item[\cite{bG06}, \cite{BS13}] There exist geometric $(n_4)$ configurations of lines if and only if $n\ge 18$ except possibly for $n\in\{19,22,23,26,37,43\}$.
\item[\cite{BP15}] There is no geometric $(19_4)$ configuration.
\item[\cite{BP16}] There exist geometric $(37_4)$ and $(43_4)$ configurations.
\end{enumerate}
Thus for the existence of geometric $(n_4)$ configurations, only the cases $n\in \{22,23,26\}$ were open.
Using the above algorithm we can produce examples when $n$ is $22$ and $26$.

We will denote both projective lines and points with coordinates $(a:b:c)$ since points and lines are dual to each other in the plane.

\subsection{$(22_4)$}

The key idea to obtain $(n_4)$ configurations with the above algorithm is to choose an arrangement $\Ac_0$ which already has some points of multiplicity $4$. This way, the orbit $\Ac$ is likely to have a large number of quadruple points as well.
Indeed, starting with an arrangement in $\FF_{19}\PP^2$ with two quadruple points and a group $H$ of order $4$, we find the following arrangement of lines (see Figure \ref{fig_22_4} and \ref{fig_22_4_e}):
\begin{eqnarray*}
\Ac_{22_4}&=&\{
(1:0:0),(0:1:0),(0:0:1),(1:1:1),(24:-5 w - 13:0), \\
&& (24:5 w + 13:24 w),(1:0:w),(2:0:w),\\
&& (24:-5 w - 13:-4 w + 52),(24:5 w + 13:28 w - 52),\\
&& (6:-w + 13:-w + 13),(24:-5 w - 13:16 w + 104),\\
&& (48:w + 65:24 w),(24:5 w + 13:-32 w + 104),\\
&& (18:-w + 13:4 w + 26),(12:-w + 13:0),\\
&& (96:w + 65:56 w - 104),(48:w + 65:-8 w + 104),\\
&& (48:w + 65:20 w + 52),(39:-w + 52:-w + 52),\\
&& (4:w + 13:4 w),(24:w + 26:12 w)
\}
\end{eqnarray*}
where $w$ is a root of $x^2+7x-26$. Each of the $22$ lines has $13$ intersection points, $4$ quadruple and $9$ double points.
The dual configuration (in which the $22$ quadruple points are the lines) has $12$ lines with $4$ quadruple, one triple, and $7$ double points, and $10$ lines with $4$ quadruple and $9$ double points (see Figure \ref{fig_22_4} and \ref{fig_22_4_e}).

\begin{remar}
\begin{enumerate}
\item Since there are two roots $w$ of $x^2+7x-26$, we obtain two arrangements $\Ac_{22_4}$ up to projectivities. The corresponding matroids are isomorphic, but the CW complexes are different. This is why we find four arrangements including the duals.
\item The corresponding matroid has a group of symmetries isomorphic to $\ZZ/2\ZZ\times\ZZ/2\ZZ$. This rather small group is probably the reason why this example did not appear in an earlier publication.
\item The above search finds these examples within a few seconds. The difficulty in finding such a configuration with the above algorithm is thus not about optimizing code.
\end{enumerate}
\end{remar}

\subsection{$(26_4)$}

The same technique yields the following $(26_4)$ configuration (and its dual), see Figure \ref{fig_26_4}:
\begin{tiny}
\begin{eqnarray*}
\Ac_{26_4}&=&\{
(1:0:0),(0:1:0),(0:0:1),(1:1:1), \\
&& (1:-z^2 - 2z:z),(z:-z^2 - 2z:z),(-1:-z:-z), \\
&& (2z^2 + 2z:-6z^2 - 14z:-5z^2 - 4z + 21), \\
&& (10z^2 + 12z - 14:-6z^2 - 56z - 98:16z^2 + 24z - 28), \\
&& (6z^2 - 14:-2z^2 - 12z + 14:-16z^2 - 20z + 56), \\
&& (68z^2 + 56z - 196:20z^2 + 72z - 84:84z^2 + 100z - 280), \\
&& (-24z^2 - 60z - 28:0:-26z^2 - 40z + 14), \\
&& (-256z^2 - 112z + 784:-352z^2 - 624z + 336:-264z^2 - 224z + 392), \\
&& (0:-16z^2 + 4z + 28:-2z^2 + 12z - 14), \\
&& (68z^2 + 56z - 196:20z^2 + 72z - 84:20z^2 + 72z - 84), \\
&& (-1136z^2 - 256z + 3696:-3152z^2 - 2560z + 7952:-1840z^2 - 528z + 4928), \\
&& (-608z^2 - 1760z + 1792:1120z^2 - 1824z - 8064:-2624z^2 - 4064z + 6048), \\
&& (0:-1120z^2 + 1824z + 8064:-1872z^2 - 1120z + 7056), \\
&& (-12864z^2 - 14976z + 48832:-23616z^2 - 35968z + 90048:-27584z^2 - 19200z + 88256), \\
&& (4288z^2 + 37888z + 61376:-44736z^2 - 170752z - 157248:-2656z^2 + 19712z + 48608), \\
&& (8z^2 + 136z - 224:-304z^2 - 392z + 1176:-412z^2 - 632z + 1652), \\
&& (-784z^2 - 608z + 2800:-272z^2 - 288z - 1232:-1136z^2 - 256z + 3696), \\
&& (11264z^2 + 30464z - 75264:-193536z^2 - 190208z + 637952:-123776z^2 - 57344z + 307328), \\
&& (-65984z^2 - 13056z + 231616:-55360z^2 - 91904z + 37184:448z^2 + 55808z - 20160), \\
&& (8192z^2 - 31232z - 155904:55808z^2 + 147712z - 57344:-54912z^2 - 36096z + 17024), \\
&& (627968z^2 + 367104z - 1732864:2495232z^2 + 3188224z - 9105152:\\
&& 1453568z^2 + 1928704z - 5465600)
\}
\end{eqnarray*}
\end{tiny}
where $z$ is the real root of $x^3 + 3x^2 - x - 7$.

\begin{remar}
All the matroids presented in this note have realizations which are unique up to projectivities and Galois automorphisms.
For $\Ac_{26_4}$ there is a complex realization which may not be transformed into a real arrangement by a projectivity, namely when $z$ is a complex root of $x^3 + 3x^2 - x - 7$.
\end{remar}

\subsection{$(23_4)$}

The arrangement of lines
\begin{eqnarray*}
\Ac_{23_4}&=&\{
(0:0:1),(0:1:0),(1:0:0),(2:0:1),(1:0:1), \\
&& (1:-1:1),(1:1:1),(2:2:i + 1),(1:1:i),(1:-i:0), \\
&& (2:-2i:i + 1),(1:-i:i + 1),(1:-i + 2:i), \\
&& (5:-3i + 4:i + 2),(2:-i + 1:i + 1),(5:-2i + 1:i + 2), \\
&& (5:-i - 2:i + 2),(5:-i + 2:-i + 2),(5:-i + 2:i + 3), \\
&& (5:-i + 2:3i + 4),(1:i:0),(1:i:-i),(1:i:i)
\}
\end{eqnarray*}
where $i=\sqrt{-1}$ has $25$ intersection points of multiplicity $4$. The right choice of $23$ points yields a $(23_4)$ configuration in the complex projective plane.

\begin{remar}
Notice that the above algorithm produces many more non isomorphic examples over finite fields and even (at least) three more examples over the complex numbers. Thus these results give no hint concerning the existence of geometric $(23_4)$ configurations.
\end{remar}

\bigskip
\textbf{Acknowledgement.}
I would like to thank J.~Bokowski and V.~Pilaud for calling my attention to the subject of $(n_k)$ configurations.

\def\cprime{$'$}
\providecommand{\bysame}{\leavevmode\hbox to3em{\hrulefill}\thinspace}
\providecommand{\MR}{\relax\ifhmode\unskip\space\fi MR }
\providecommand{\MRhref}[2]{%
  \href{http://www.ams.org/mathscinet-getitem?mr=#1}{#2}
}
\providecommand{\href}[2]{#2}


\begin{thebibliography}{1}

\bibitem{BP15}
J\"urgen Bokowski and Vincent Pilaud, \emph{On topological and geometric
  {$(19_4)$} configurations}, European J. Combin. \textbf{50} (2015), 4--17.

\bibitem{BP16}
\bysame, \emph{Quasi-configurations: building blocks for point-line
  configurations}, Ars Math. Contemp. \textbf{10} (2016), no.~1, 99--112.

\bibitem{BS13}
J\"urgen Bokowski and Lars Schewe, \emph{On the finite set of missing geometric
  configurations {$(n_4)$}}, Comput. Geom. \textbf{46} (2013), no.~5, 532--540.

\bibitem{p-C10b}
M.~Cuntz, \emph{Minimal fields of definition for simplicial arrangements in the
  real projective plane}, Innov. Incidence Geom. \textbf{12} (2011), 49--60.

\bibitem{bG06}
Branko Gr\"unbaum, \emph{Connected {$(n_4)$} configurations exist for almost
  all {$n$}---second update}, Geombinatorics \textbf{16} (2006), no.~2,
  254--261.

\end{thebibliography}
\end{document}